\documentclass[a4paper,11pt]{amsart}
\usepackage{amsmath,amssymb,amsfonts,latexsym}
\usepackage[top=5cm, bottom=5cm, left=4cm, right=4cm]{geometry}
\usepackage[colorlinks=true,urlcolor=magenta,linkcolor=red,citecolor=blue]{hyperref}

\newtheorem{theorem}{Theorem}[section]

\newtheorem{proposition}[theorem]{Proposition}

\theoremstyle{remark}
\newtheorem{remark}[theorem]{Remark}

\input{tcilatex}

\begin{document}
\title[On the Legendre polynomials]{The Legendre
polynomials associated with Bernoulli, Euler, Hermite and Bernstein polynomials}
\author[S. Araci]{Serkan Araci}
\address{Atat\"{u}rk Street, 31290 Hatay, TURKEY}
\email{mtsrkn@hotmail.com}
\author[M. Acikgoz]{Mehmet Acikgoz}
\address{University of Gaziantep, Faculty of Science and Arts, \\
Department of Mathematics, 27310 Gaziantep, TURKEY}
\email{acikgoz@gantep.edu.tr}
\author[A. Bagdasaryan]{Armen Bagdasaryan}
\address{Russian Academy of Sciences, Institute for Control
Sciences, 65 Profsoyuznaya, 117997 Moscow, RUSSIA}
\email{abagdasari@hotmail.com}
\author[E. \c{S}en]{Erdo\u{g}an \c{S}en}
\address{Department of Mathematics, Faculty of Science and Letters,
Namik Kemal University, 59030 Tekirda\u{g}, TURKEY}
\email{erdogan.math@gmail.com}

\begin{abstract}
In the present paper, we deal mainly with arithmetic properties of Legendre
polynomials by using their orthogonality property. We
show that Legendre polynomials are proportional with Bernoulli, Euler,
Hermite and Bernstein polynomials.

\vspace{2mm}

\noindent \textsc{Keywords and phrases.} Legendre polynomials, Bernoulli
polynomials, Euler polynomials, Hermite \ polynomials, Bernstein
polynomials, orthogonality.
\end{abstract}

\thanks{{Published in: {\em Turkish J Analysis Number Theory} {\bf 1} (2013) 1--3.}}

\maketitle




\section{Introduction}


Legendre polynomials, which are special cases of Legendre functions, are
introduced in 1784 by the French mathematician A. M. Legendre (1752-1833).
Legendre functions are a vital and important in problems including spherical
coordinates. Due to their orthogonality properties they are also useful in
numerical analysis (see \cite{Andrews}). Besides, the Legendre polynomials, $%
P_{n}\left( x\right) $, are described via the following generating function:%
\begin{equation}
\frac{1}{\sqrt{1-2xt+t^{2}}}=\sum_{n=0}^{\infty }P_{n}\left( x\right) t^{n}%
\text{.}  \label{equation 1}
\end{equation}

Legendre polynomials are the everywhere regular solutions of \textit{Legendre's differential equation} that we can write as follows:%
\begin{multline*}
\left(1-x^2\right) \frac{d}{dx}P_n(x) - 2x\frac{d}{dx}P_n(x) + m P_n(x) = \\
                             = \frac{d}{dx}\left[\left(1-x^2\right) \frac{d}{dx}P_n(x)\right] + m P_n(x) = 0,
\end{multline*}
where $m=n\left( n+1\right) $ and $n=0,1,2,\cdots .$ Taking $x=1$ in (\ref%
{equation 1}) and by using geometric series, we see that $P_{n}\left(
1\right) =1$, so that the Legendre polynomials are $normalized$.

Legendre polynomials can be generated using \textit{Rodrigue's formula} as
follows:%
\begin{equation}
P_{n}\left( x\right) =\frac{1}{n!2^{n}}\frac{d^{n}}{dx^{n}}\left(
x^{2}-1\right) ^{n}\text{.}  \label{equation 2}
\end{equation}

Note that the right hand side of (\ref{equation 2}) is a polynomial (see 
\cite{Bailey}, \cite{Andrews}).

The Bernoulli polynomials are defined by means of the following generating
function: 
\begin{equation}
\sum_{n=0}^{\infty }B_{n}\left( x\right) \frac{t^{n}}{n!}=\frac{t}{e^{t}-1}%
e^{xt},\text{ }\left\vert t\right\vert <2\pi \text{ (see \cite{Kellner}, 
\cite{Kim 6}).}  \label{equation 3}
\end{equation}

By (\ref{equation 3}), we know that $\frac{dB_{n}\left( x\right) }{dx}%
=nB_{n-1}\left( x\right) $. Taking $x=0$ in (\ref{equation 3}), we have $%
B_{n}\left( 0\right) :=B_{n}$ that stands for $n$-th Bernoulli number. 

The Euler polynomials are known to be defined as:%
\begin{equation}
\sum_{n=0}^{\infty }E_{n}\left( x\right) \frac{t^{n}}{n!}=\frac{2}{e^{t}+1}e^{xt}.  \label{equation 4}
\end{equation}


The Euler polynomials can also be expressed by explicit formulas, e.g. 
\begin{equation*}
E_{n}(x)=\sum_{k=0}^{n}\binom{n}{k}\frac{E_{k}}{2^{k}}\left( x-\frac{1}{2}\right) ^{n-k},
\end{equation*}%
where $E_{k}$ means the Euler numbers. These numbers are
expressed with the Euler polynomials through $E_{k}=2^{k}E_{k}(1/2)$. 

Now also, we give the definition of Hermite polynomials as follows: 
\begin{equation}
e^{2xt-t^{2}}=\sum_{n=0}^{\infty }H_{n}\left( x\right) \frac{t^{n}}{n!}\text{%
.}  \label{equation 5}
\end{equation}

Let $C\left( \left[ 0,1\right] \right) $ be the space of continuous
functions on $\left[ 0,1\right] $. For $f\in C\left( \left[ 0,1\right]
\right) $, Bernstein operator for $f$ is defined by 
\begin{equation*}
\mathcal{B}_{n}\left( f,x\right) =\sum_{k=0}^{n}f\left( \frac{k}{n}\right)
B_{k,n}\left( x\right) =\sum_{k=0}^{n}f\left( \frac{k}{n}\right) \binom{n}{k}%
x^{k}\left( 1-x\right) ^{n-k}\text{,}
\end{equation*}%
where $n,$ $k\in 
\mathbb{N}
^{\ast }:=%
\mathbb{N}
\cup \left\{ 0\right\} $ and $%
\mathbb{N}
$ is the set of natural numbers. Here $B_{k,n}\left( x\right) $ is called
Bernstein polynomials, which are defined by 
\begin{equation}
B_{k,n}\left( x\right) =\binom{n}{k}x^{k}\left( 1-x\right) ^{n-k},\text{ }%
x\in \left[ 0,1\right] \text{ (cf. \cite{Araci 3}, \cite{Kim7}.)}
\label{equation 6}
\end{equation}

In \cite{Andrews}, \cite{Bailey}, the orthogonality of Legendre polynomials
is known as 
\begin{equation}
\int_{-1}^{1}P_{m}\left( x\right) P_{n}\left( x\right) dx=\frac{2}{2n+1}%
\delta _{m,n},\text{where }\delta _{m,n}\text{ is Kronecker's delta.}
\label{equation 7}
\end{equation}

In \cite{Kim 6}, by using orthogonality property of Legendre \cite{Kim 6},
Kim $et$ $al.$ effected interesting identities for them. We also obtain some
interesting properties of the Legendre polynomials arising from Bernoulli,
Euler, Hermite and Bernstein polynomials.

\section{Identities on the Legendre polynomials arising from
Bernoulli, Euler, Hermite and Bernstein polynomials}

Let $\mathcal{P}_{n}=\left\{ q\left( x\right) \in 
\mathbb{Q}
\left[ x\right] \mid \deg p\left( x\right) \leq n\right\} $. Then we define
an inner product on $\mathcal{P}_{n}$ as follows:%
\begin{equation}
\left\langle q_{1}\left( x\right) ,q_{2}\left( x\right) \right\rangle
=\int_{-1}^{1}q_{1}\left( x\right) q_{2}\left( x\right) dx\text{, }\left(
q_{1}\left( x\right) ,q_{2}\left( x\right) \in \mathcal{P}_{n}\right) \text{.%
}  \label{equation 9}
\end{equation}

Note that $P_{0}\left( x\right) ,P_{1}\left( x\right) ,...,P_{n}\left(
x\right) $ are the orthogonal basis for $\mathcal{P}_{n}$. Let us now
consider $q\left( x\right) \in \mathcal{P}_{n}$; then we see that%
\begin{equation}
q\left( x\right) =\sum_{k=0}^{n}C_{k}P_{k}\left( x\right),
\label{equation 10}
\end{equation}
where the coefficients $C_k$ are defined over the field of real numbers.

From the above, we readily see that%
\begin{eqnarray}
C_{k} &=&\frac{2k+1}{2}\left\langle q\left( x\right) ,P_{k}\left( x\right)
\right\rangle =\frac{2k+1}{2}\int_{-1}^{1}P_{k}\left( x\right) q\left(
x\right) dx  \label{equation 11} \\
&=&\frac{2k+1}{k!2^{k+1}}\int_{-1}^{1}\left( \frac{d^{k}}{dx^{k}}\left(
x^{2}-1\right) ^{k}\right) q\left( x\right) dx\text{.}  \notag
\end{eqnarray}

By (\ref{equation 10}) and (\ref{equation 11}), we have the following
proposition.

\begin{proposition}
\label{propstion1}Let $q\left( x\right) \in \mathcal{P}_{n}$ and $q\left(
x\right) =\sum_{k=0}^{n}C_{k}P_{k}\left( x\right) $, then 
\begin{equation*}
C_{k}=\frac{2k+1}{k!2^{k+1}}\int_{-1}^{1}\left( \frac{d^{k}}{dx^{k}}\left(
x^{2}-1\right) ^{k}\right) q\left( x\right) dx \text{ (see \cite{Kim 6}).}
\end{equation*}
\end{proposition}

If we take $q\left( x\right) =x^{n}$ in Proposition (\ref{propstion1}), the
coefficients $C_{k}$ can be found as%
\begin{equation}
C_{k}=\frac{\left( 2k+1\right) 2^{k+1}}{\left( n+k+2\right) !}\frac{n!\left( 
\frac{n+k+2}{2}\right) !}{\left( \frac{n-k}{2}\right) !}\text{ for }%
n-k\equiv 0\left( \func{mod}2\right) \text{ (see \cite{Kim 6}).}
\label{equation 15}
\end{equation}

Let $q\left( x\right) =B_{n}\left( x\right) $. Then by using Proposition \ref%
{propstion1} and (\ref{equation 15}), we have%
\begin{equation*}
C_{k}=\frac{2k+1}{k!2^{k+1}}\int_{-1}^{1}\left( \frac{d^{k}}{dx^{k}}\left(
x^{2}-1\right) ^{k}\right) B_{n}\left( x\right) dx
\end{equation*}%
where $B_{n}\left( x\right) $ are the aforementioned Bernoulli polynomials that
can be expressed through Bernoulli numbers $B_{n}$ as follows:%
\begin{equation*}
B_{n}\left( x\right) =\sum_{j=0}^{n}\binom{n}{j}B_{n-j}x^{j}.
\end{equation*}%
From this, we have%
\begin{multline*}
C_{k} =\sum_{j=0}^{n}\binom{n}{j}B_{n-j}\left[ \int_{-1}^{1}\left( \frac{%
d^{k}}{dx^{k}}\left( x^{2}-1\right) ^{k}\right) x^{j}dx\right]  \\
=\left( 2^{k+2}k+2^{k+1}\right) \sum_{j=0}^{n}\frac{j!\binom{n}{j}\left( 
\frac{j+k+2}{2}\right) !}{\left( \frac{j-k}{2}\right) !\left( j+k+2\right) !}%
B_{n-j}\text{ for }j-k\equiv 0\left( \func{mod}2\right) \text{.}
\end{multline*}

Therefore we have the following theorem.

\begin{theorem}
Let $B_{n}\left( x\right) =\sum_{k=0}^{n}C_{k}P_{k}\left( x\right) \in 
\mathcal{P}_{n}$. Then we have 
\begin{multline*}
B_{n}\left( x\right) =2\sum_{k=0}^{n}\Biggl( \left( 2^{k+2}k+2^{k+1}\right) \times \\
        \times \sum_{j-k\equiv 0\left( \func{mod}2\right) }^{n}\frac{j!\binom{n}{j}\left( 
         \frac{j+k+2}{2}\right) !}{\left( \frac{j-k}{2}\right) !\left( j+k+2\right) !}%
          B_{n-j}\Biggr) P_{k}\left( x\right) \text{.}
\end{multline*}
\end{theorem}

Let $H_{n}\left( x\right) \in \mathcal{P}_{n}$. By Proposition \ref%
{propstion1} and (\ref{equation 15}), we have the following theorem.

\begin{theorem}
Let $H_{n}\left( x\right) =\sum_{k=0}^{n}C_{k}P_{k}\left( x\right) \in 
\mathcal{P}_{n}$. Then we have%
\begin{multline*}
H_{n}\left( x\right) =\sum_{k=0}^{n}\Biggl( \left( 2^{k+2}k+2^{k+1}\right)  \times \\
       \times   \sum_{j-k\equiv 0\left( \func{mod}2\right) }^{n}\frac{2^{j}\binom{n}{j}%
                j!\left( \frac{j+k+2}{2}\right) !}{\left( j+k+2\right) !\left( \frac{j-k}{2}%
                \right) !}H_{n-j}\text{ }\Biggr) P_{k}\left( x\right) \text{.}
\end{multline*}
\end{theorem}

Let the Bernstein polynomials $B_{j,n}\left( x\right) \in \mathcal{P}_{n}$. By Proposition \ref%
{propstion1} and (\ref{equation 15}), we have the following theorem.

\begin{theorem}
Let $B_{j,n}\left( x\right) =\sum_{k=0}^{n}C_{k}P_{k}\left( x\right) \in 
\mathcal{P}_{n}$.  We have%
\begin{multline*}
B_{j,n}\left( x\right) =\sum_{k=0}^{n}\Biggl( \left( 2^{k+2}k+2^{k+1}\right) \times \\
        \times  \sum_{l+j-k\equiv 0\left( \func{mod}2\right) }^{n-j}\frac{\binom{n-j}{l}%
                \left(-1\right)^{l}\left( l+j\right) !\left( \frac{l+j+k+2}{2}\right) !}{%
                \left(l+j+k+2\right) !\left( \frac{l+j-k}{2}\right) !}\Biggr) P_{k}\left(x\right).
\end{multline*}
\end{theorem}

The following equality is defined by Kim et al. in \cite{Kim 6}:%
\begin{equation}
\sum_{k=0}^{n}B_{k}\left( x\right) B_{n-k}\left( x\right) =\frac{2}{n+2}%
\sum_{l=0}^{n-2}\binom{n+2}{l}B_{n-l}B_{l}\left( x\right) +\left( n+1\right)
B_{n}\left( x\right) \text{.}  \label{equation 14}
\end{equation}

Let $\sum_{k=0}^{n}B_{k}\left( x\right) B_{n-k}\left( x\right) \in \mathcal{P%
}_{n}$. By Proposition \ref{propstion1} and (\ref{equation 15}), we get the
following theorem.

\begin{theorem}
Let $\sum_{k=0}^{n}B_{k}\left( x\right) B_{n-k}\left( x\right) \in \mathcal{P%
}_{n}$. Then we have%
\begin{multline*}
\sum_{k=0}^{n}B_{k}\left( x\right) B_{n-k}\left( x\right) 
=\sum_{k=0}^{n}\left( 2^{k+2}k+2^{k+1}\right) \times \\ 
\times \Biggl[\frac{2}{n+2}\sum_{l=0}^{n-2}\sum_{j-k\equiv 0\left( \func{mod}2\right)
}^{l}B_{n-l}B_{l-j} \frac{\binom{n+2}{l}\binom{l}{j}j!\left( \frac{j+k+2}{2}\right) !}{%
\left( \frac{j-k}{2}\right) !\left( j+k+2\right) !} \\
+\left( n+1\right) \sum_{l-k\equiv 0\left( \func{mod}2\right) }^{n}\binom{n%
}{l}B_{n-l}\frac{l!\left( \frac{l+k+2}{2}\right) !}{\left( \frac{l-k}{2}%
\right) !\left( l+k+2\right) !}\Biggr] P_{k}\left( x\right) \text{.}
\end{multline*}
\end{theorem}

Let $q\left( x\right) =\sum_{k=0}^{n}E_{k}\left( x\right) E_{n-k}\left(
x\right) \in \mathcal{P}_{n}$. In \cite{Kim8}, Kim et al derived convolution
formula for the Euler polynomials as 
\begin{equation*}
\sum_{k=0}^{n}E_{k}\left( x\right) E_{n-k}\left( x\right) =-\frac{4}{n+2}%
\sum_{l=0}^{n}\binom{n+2}{l}E_{n-l+1}B_{l}\left( x\right) \text{.}
\end{equation*}

By Proposition \ref{propstion1} and (\ref{equation 15}), we get the
following theorem.

\begin{theorem}
The following equality holds true:%
\begin{eqnarray*}
\sum_{k=0}^{n}E_{k}\left( x\right) E_{n-k}\left( x\right)  &=&-\frac{8}{n+2}%
\sum_{k=0}^{n}\left( 2^{k+1}k+2^{k}\right)  \\
&&\times \{\sum_{l=0}^{n}\sum_{j-k\equiv 0\left( \func{mod}2\right) }^{l}%
\binom{n+2}{l}\binom{l}{j}E_{n-l+1} \\
&&\times B_{l-j}\frac{j!\left( \frac{j+k+2}{2}\right) !}{\left( \frac{j-k}{2}%
\right) !\left( j+k+2\right) !}\}P_{k}\left( x\right) \text{.}
\end{eqnarray*}
\end{theorem}

\begin{remark}
By using Theorem \ref{propstion1}, we can find many interesting identities
for the special polynomials in connection with Legendre polynomials.
\end{remark}

%


\end{document}